\input amstex
\documentstyle{amsppt}
\loadmsbm
\hcorrection{18mm}  \vcorrection{5mm}
\nologo \TagsOnRight \NoBlackBoxes
\headline={\ifnum\pageno=1 \hfill\else%
{\tenrm\ifodd\pageno\rightheadline \else
\leftheadline\fi}\fi}
\def\rightheadline{EJDE--2001/??\hfil Two symmetry problems in potential theor
\hfil\folio}
\def\leftheadline{\folio\hfil Tewodros Amdeberhan
 \hfil EJDE--2001/??}

\topmatter
\title
Two symmetry problems in potential theory
\endtitle

\thanks 
\hfil\break\indent
\copyright 2001 Southwest Texas State University. \hfil\break\indent
Submitted February 13, 2001. Published ??.
\endthanks

\author  Tewodros Amdeberhan  \endauthor
\address DeVry College of Technology, 
North Brunswick NJ 08902, USA 
\endaddress
\email amdberha\@nj.devry.edu, tewodros\@math.temple.edu
\endemail

\abstract
We consider two overdetermined boundary value problems
as variants on J. Serrin's 1971 classical results and prove
in both cases that the domains must be Euclidean balls.
\endabstract
\endtopmatter

\document

Assume throughout that $\Omega \subseteq {\Bbb R}^N$ is a
bounded domain whose boundary is of class $C^2$ and contains the
origin strictly in its  interior.  Let $\nu$ be the outer
unit normal to $\partial \Omega.$ Summation over repeated indices
is in effect. 
In 1971, James Serrin [1] proved the following classical result.

\proclaim{Theorem 1} Suppose there exists a function $u \in C^2(\bar\Omega)$ 
satisfying the elliptic differential equation 
$$a(u,|p|)\Delta u + h(u,|p|)u_iu_ju_{ij}=f(u,|p|) \quad \text{in $\Omega$}
$$ 
where $a, f$ and $h, p_i, p_j$ are continuously differentiable functions 
of $u$ and $p$ (here $p=(u_1,\dots,u_n)$ denotes the gradient vector of $u$).
Suppose also that $u >0$ in $\Omega$ and that $u$ satisfies the boundary 
conditions 
$$u=0, \quad \frac{\partial u}{\partial \nu}=\text{constant} \quad 
\text{on $\partial \Omega$.}
$$
Then $\Omega$ must be a ball and $u$ is radially symmetric.  
\endproclaim

The method of proof combines the Maximum Principles and the device 
(which goes back to A.  D.  Alexandroff:  {\it every embedded surface in 
${\Bbb R}^N$ with constant mean curvature must be a sphere}) 
of {\it moving planes} to a critical position and then showing that the 
solution is symmetric about the limiting plane.  
In a subsequent article, H.  F.  Weinberger [3] gave a
simplified proof for the special case of the Poisson differential equation,
$\Delta u = -1$.

Our aim at present is to introduce some variants on Serrin's result and arrive
at the same symmetry conclusions by employing elementary arguments.  The next
statement involves radial dependence on the boundary conditions.

\proclaim{Proposition 1} 
Suppose there exists a solution $u \in C^2(\bar\Omega)$ to the overdetermined 
problem:
$$ \gathered
\Delta u = -1 \quad \text{in $\Omega$} \\
u = 0 \quad \text{on $\partial \Omega$} \\
\frac {\partial u}{\partial \nu} = -cr 
\quad \text{on $\partial \Omega$;} \endgathered \tag1 
$$
where $r=\sqrt {x_1^2+...+x_N^2}$ and $c$ is a constant. 
Then $\Omega$ is an $N$-dimensional ball.
\endproclaim

Before turning to the details we like to discuss some physical motivations for
the problem itself.  Consider a viscous incompressible fluid moving in straight
parallel streamlines through a straight pipe of given cross-sectional form
$\Omega$.  If we fix rectangular coordinates in space with the $z-$axis directed
along the pipe, it is well known that the flow velocity $u$ is then a function
of $x, y$ alone satisfying the Poisson differential equation (for $N=2$)
$\Delta u = - A$ in $\Omega$, where $A$ is a constant related to the 
viscosity of the fluid and to the rate of change of pressure per unit 
length along the pipe.
Supplementary to the differential equation one has the adherence condition $u=0$
on $\partial \Omega$.

Finally, the tangential stress per unit area on the pipe wall is given by the 
quantity $\mu \frac{\partial u}{\partial \nu}$, where $\mu$ is the viscosity. 
Our proposition then states that  the ratio of the tangential stress on the 
pipe wall to its radial distance is the same
at all points of the wall if and only if the pipe has a circular cross section. \rm

Exactly the same differential equation and boundary condition arise in the
linear theory of torsion of a solid straight bar of cross-section $\Omega$; see
[2] pp.109-119.  In light of this, Proposition 1 states that,  when a solid
straight bar is subject to torsion, the ratio of the magnitude of the resulting
traction, which occurs at the surface of the bar, to the radial distance of the
surface is independent of position if and only if the bar has a circular
cross-section.

\proclaim{Lemma 1}
Under the hypothesis of Proposition 1, the following holds
$$
\int_{\Omega}u dx = c^2\int_{\Omega}r^2 dx.\tag2
$$
\endproclaim

\demo{Proof}
Let us introduce the auxiliary function $h = 2u - x_iu_i$. Then clearly
$\Delta h =0$. Green's identities
$$
\int_{\Omega}\left(h\Delta u - u\Delta h\right) dx
=\int_{\partial \Omega}\left(h\frac {\partial u}{\partial \nu} -
u\frac {\partial h}{\partial \nu}\right) d\sigma, 
$$
The boundary conditions on $u$, and the harmonicity of $h$ lead to
$$
\int_{\Omega} h dx = c\int_{\partial \Omega} rh d\sigma \tag3 
$$
Now we compute the left and right hand sides of (3) individually.
Applying the divergence theorem, we get
$$
\aligned
\int_{\Omega} h dx =& \int_{\Omega}(2+N)u dx - \int_{\Omega} div(xu) dx \\
=&\int_{\Omega} (2+N)u dx - \int_{\partial \Omega} u x\cdot \nu d\sigma \\
=&\int_{\Omega} (2+N)u dx, \quad\text{for $u=0$ on $\partial \Omega.$} 
\endaligned \tag4
$$
Again since $u$ vanishes on the boundary $\partial\Omega$
(therefore $\nu=\pm \frac{\nabla u}{\Vert \nabla u \Vert}$) and 
$r_i=\frac{x_i}r$, we gather that
(note: the argument here could have proceeded using the so-called 
Pohozaev's identity, but we drop it so as not to use any ``heavy gun'')
$$
\aligned
c\int_{\partial \Omega} rh d\sigma 
=& -c\int_{\partial \Omega} rx_iu_i d\sigma 
 = -c\int_{\partial \Omega} r^2\frac {\partial u}{\partial r}d\sigma\\
=&-c\int_{\partial \Omega} r^2 \frac {\partial u}{\partial \nu}
 \frac {\partial r}{\partial \nu} d\sigma 
 =c^2\int_{\partial \Omega} r^3 \frac {\partial r}{\partial \nu} d\sigma \\
=&\frac {c^2}4\int_{\partial \Omega} \frac {\partial (r^4)}{\partial \nu} 
  d\sigma
=\frac {c^2}4\int_{\Omega} \Delta (r^4) dx \\
=&c^2(N+2)\int_{\Omega} r^2 dx. 
\endaligned \tag5
$$
Combining (3), (4) and (5) proves the Lemma. 
\qed\enddemo

\demo{Proof of Proposition 1}
Consider the functional $\varPhi = u_iu_i - c^2r^2$. Then,
$$
\aligned
\Delta \varPhi &= 2u_{ij}u_{ij} + 2u_i\Delta u_i - 2c^2N \\
&= 2u_{ij}u_{ij} - 2c^2N\quad \text{since $\Delta$u is a constant}\\
&=2\sum_{i,j}\left(u_{ij}+\frac {\delta_{ij}}N\right)^2 +\frac 2N - 2c^2N\\
&\geq \frac 2N - 2c^2N. \endaligned\tag6
$$
For a moment assume that $cN \leq 1$.  Then we have,
$$\Delta \varPhi\geq 0 \quad\text{in }\Omega.\tag7
$$
Note that $\nu =\pm \nabla u/ \Vert \nabla u\Vert$, since $u$ vanishes
on the boundary. Thus on $\partial \Omega$ we have
$$\varPhi = \left(\frac {\partial u}{\partial \nu}\right)^2- c^2r^2 =
0.\tag8$$
Applying Green's identities and replacing the boundary conditions results in
$$
\aligned
\int_{\Omega} \varPhi dx &= \int_{\Omega} \nabla u \cdot \nabla u -c^2\int_{\Omega} r^2 dx \\
&=\int_{\partial \Omega} u\frac {\partial u}{\partial \nu} d\sigma - \int_{\Omega} u\Delta u dx-c^2\int_{\Omega} r^2 dx \\
&=\int_{\Omega} u dx - c^2\int_{\Omega} r^2 dx. \endaligned \tag9 
$$
By Lemma 1 above, equation (9) yields
$$\int_{\Omega} \varPhi = 0.$$
Standard Maximum Principles together with the properties (7)-(9) of 
$\varPhi$ imply that $\varPhi \equiv$  0 in $\Omega$. This in particular 
forces equality in (6), i.e.
$$u_{ij} + \frac {\delta_{ij}}N \equiv 0.$$
Consequently, $u$ takes the radial form
$$u = a - \frac {r^2}{2N}.$$
Again, since $u$ vanishes on the boundary $\partial \Omega$, we obtain that
$\Omega$ is indeed a ball, as asserted. 
\qed\enddemo

\subhead Remarks\endsubhead
\roster
\item  It is not hard to see that the $C^2$-smoothness assumptions, on $u$, 
made in the preceding Proposition (even Proposition 2, below) could be 
relaxed except the proofs would then get rather technical. 

\item The assumption $cN \leq 1$ in the proof of Proposition 1 is not essential.
To see this, notice in fact that due to the conditions (1) on $u$ coupled 
with the divergence  theorem verify that
$$
\aligned
|\Omega | &= - \int_{\Omega} \Delta u dx 
= - \int_{\partial \Omega} \frac {\partial u}{\partial \nu} d\sigma 
= c\int_{\partial \Omega} r d\sigma \\
&\geq cN\int_{\partial \Omega} \frac {r dA}N = cN | \Omega | \quad 
\text{from geometry.}
\endaligned
$$
Therefore, $cN \leq $ 1.

\item Using the same argument and replacing the second boundary condition by 
$\frac {\partial u}{\partial \nu} = -cr^{\alpha}$, where $\alpha\geq1$ 
(note: this already implies that 
$cN\operatorname{dist}(0,\partial\Omega)^{\alpha-1} \leq 1$), 
$d=\operatorname{diam}(\Omega)$ and
$$0 < c \leq \frac {2^{\alpha-1}}{ d^{\alpha-1}\sqrt{\alpha N(N+2\alpha-2)}}\*
$$
we still get the conclusion of Proposition 1.

\item Even more generally, if $\frac {\partial u}{\partial \nu} = -f(r)$ where
$g(r)=f^2(r)$ satisfies the conditions 
$$\gathered
\int_{\Omega}(rg'-2g)dx \geq 0, \\
g''+\frac {N-1}r g' - \frac 2N \leq 0
 \endgathered$$
then  Proposition 1  holds. Moreover, it turns out that either 
$f\equiv c_1$ or $f = c_2r$ for some positive constants $c_1$ and $c_2$. 
One of the implications of which is that under these general suppostions, 
there can only be two possible forms for the boundary derivative: either 
$\partial u/\partial \nu =-c_1$ as in Theorem 1 of Serrin, or 
$\partial u/\partial \nu =-c_2r$ as in Proposition 1
of this article!
\endroster

\proclaim{Proposition 2}
Suppose there exists a solution $u \in C^2(\bar\Omega)$ to following 
overdetermined boundary-value problem
$$\gathered
\Delta u = -r^\alpha \quad \text{in $\Omega$}\\
u = 0 \quad \text{on $\partial \Omega$} \\
x_iu_i + c_0r^{2+\alpha} + c_1 = 0 \quad \text{on $\partial \Omega$;}
\endgathered\tag10
$$
where $\alpha$, $c_0$ and $c_1$ are constants. If 
$\beta := (\alpha + 2)(c_0(\alpha + N) -1)$ is  not a negative integer,
then $\Omega$ is an $N$-dimensional ball.
\endproclaim

\demo{Proof}
Introduce the functional $V = x_iu_i + c_0r^{\alpha + 2} + c_1 + \beta u$. 
Then by direct computation we obtain that V is harmonic. 
As can easily be seen, the boundary conditions
make $V \equiv 0$ on $\partial \Omega$. Thus,
$$\gathered
\Delta V = 0 \quad \text{in $\Omega$} \\
V = 0 \quad \text{on $\partial \Omega$ .}
\endgathered $$
Classical Maximum Principles show that $V \equiv 0$ inside $\Omega$, too. 
Hence, it follows that 
$$x_iu_i + c_0r^{\alpha + 2} + c_1 + \beta u = 0 \quad 
\text{in $\overline \Omega$.} \tag11
$$
Rewriting equation (11) one obtains
$$
\frac {\partial}{\partial r}\left(r^{\beta}(u-u(0)) 
+ \frac {c_0r^{\beta + \alpha +  2}}{\beta + \alpha + 2}\right) 
= r^{\beta - 1}\left(\beta (u-u(0)) + r\frac {\partial u}{\partial r} 
+ c_0r^{\alpha + 2}\right) = 0\,.
$$
This in turn implies that
$$
r^{\beta}(u-u(0)) + \frac {c_0r^{\beta + \alpha +  2}}{\beta + \alpha + 2}
\equiv G(\Theta)
$$
for angular variables $\Theta$. Hence,
$$u - u(0) = - c_0\frac {r^{\alpha + 2}}{\beta + \alpha + 2} 
+ r^{-\beta}G(\Theta) \tag12
$$
Now, if $\beta \geq$ 0 then $u$ cannot be regular at the origin unless 
$G(\Theta) \equiv 0.$ 

If $\beta < 0$ is not a negative integer, then once more we must have $
G(\Theta)\equiv 0$ for $u$ must satisfy $\Delta u = -r^{\alpha}$ and 
$r^{-\beta}G(\Theta)$ is not harmonic with this value of $\beta$.
Therefore, the solution $u$ is radial and
$$ u = u(0) - c_0\frac {r^{\alpha + 2}}{\beta + \alpha + 2}. 
$$
After using (10) and the vanishing of $u$ on $\partial\Omega$, this proves
that $\Omega$ is a ball thereby completing the proof of the proposition. 
\enddemo


\Refs
%
\ref\key 1
\by J. Serrin
\paper A symmetry problem in potential theory
\jour Arch. Rat. Mech. Anal., 43
\yr 1971
\pages 304--318
\endref
%
\ref\key 2
\by I.S. Sokolinikoff
\paper Mathematical theory of elasticty
\jour New York: McGraw Hill
\yr 1956
\endref
%
\ref\key 3
\by H.F. Weinberger
\paper Remark on the preceding paper of Serrin
\jour Arch. Rat. Mech. Anal., 43
\yr 1971
\pages 319--320
\endref
\endRefs
\enddocument